\documentclass[conference]{IEEEtran}

        \usepackage{caption}
  	\usepackage[pdftex]{graphicx}
  	\graphicspath{{../pdf/}{../jpeg/}}
	\DeclareGraphicsExtensions{.pdf,.jpeg,.png}
        \usepackage{multirow}
        \usepackage{multicol}
        \usepackage{makecell}
	\usepackage[cmex10]{amsmath}
        \usepackage{CJKutf8}
	\usepackage{mathabx}
	\usepackage{algorithmic}
	\usepackage{array}
	\usepackage{mdwmath}
	\usepackage{mdwtab}
	\usepackage{eqparbox}
	\usepackage{url}
        \usepackage{pifont}

        \usepackage[ruled,linesnumbered]{algorithm2e}

\begin{document}
\begin{CJK}{UTF8}{gbsn}

\title{\LARGE Bi-Band ECoGNet for ECoG Decoding on Classification Task}

 \author{\authorblockN{Changqing JI\authorrefmark{1}, }
 \authorblockA{\authorrefmark{1}∗Graduate School of Information Sciences University, Miyagi, Japan 
 \\Email: \{ji.changqing\}@qq.com}
}

\maketitle

\begin{abstract} 
In the application of brain-computer interface (BCI), being able to accurately decode brain signals is a critical task. For the multi-class classification task of brain signal ECoG, how to improve the classification accuracy is one of the current research hotspots. ECoG acquisition uses a high-density electrode array and a high sampling frequency, which makes ECoG data have a certain high similarity and data redundancy in the temporal domain, and also unique spatial pattern in spatial domain. How to effectively extract features is both exciting and challenging. Previous work \cite{0} found that visual-related ECoG can carry visual information via frequency and spatial domain. Based on this finding, we focused on using deep learning to design frequency and spatial feature extraction modules, and proposed a Bi-Band ECoGNet model based on deep learning. The main contributions of this paper are: 1) The Bi-BCWT (Bi-Band Channel-Wise Transform) neural network module is designed to replace the time-consume method MST \cite{-1}, this module  greatly improves the model calculation and data storage efficiency, and effectively increases the training speed; 2) The Bi-BCWT module can effectively take into account the information both in low-frequency and high-frequency domain, which is more conducive to ECoG multi-classification tasks; 3) ECoG is acquired using 2D electrode array, the newly designed 2D Spatial-Temporal feature encoder can extract the 2D spatial feature better. Experiments have shown that the unique 2D spatial data structure can effectively improve classification accuracy; 3) Compared with \cite{0}, the Bi-Band ECoGNet model is smaller and has higher performance, with an accuracy increase of 1.24\%, and the model training speed is increased by 6 times, which is more suitable for BCI applications.
\end{abstract}

\IEEEoverridecommandlockouts

\IEEEpeerreviewmaketitle


\section{Introduction}
Brain science has been one of the hottest research areas, electroencephalography is one of the technologies that effectively records brain activity. According to the category of invasive and non-invasive, there are two types of methods: EEG (electroencephalography) and ECoG (electrocorticography). EEG measuring uses a group of electrodes which attached to the surface of the scalp. It has the advantages of low cost and non-invasiveness. However, the collected data is easily affected by noise, and the electrode layout is sparse in space. ECoG signal acquisition is to place the electrode array directly on the gray matter cortex of the brain, directly collect the electrical signals of the gray matter cortex. Compared with EEG, there is no attenuation from the skull, the signal has a higher signal-to-noise ratio, and the high density of electrodes gives us a high spatial resolution. This makes ECoG have a wide range of applications. For example, ECoG is used to assist neurological rehabilitation therapy \cite{1}. The cortical temperature, and cerebral thermodynamics of ECoG signals can contribute to improving the assessment of epileptic seizures \cite{2}, epileptic seizure detection \cite{3}, mental fatigue detection \cite{31} and ECoG signals are extensively used to find the focus of seizures in epilepsy \cite{4}. ECoG signals have a strong signal-to-noise ratio, which is conducive to the extraction of signal features. However, due to the large channel number, data redundancy, and strong similarity between channel signals, decoding ECoG signals is very challenging. Currently, the processing method of ECoG signals can be roughly divided into two categories: traditional methods and deep learning method. \\
Traditional methods tend to manually design feature extractor, focusing on brain signal statistical information. In multi-classification tasks, beside basic statistical features, such as signal average, peak value, STD, etc., there are also feature extraction methods that target at temporal-frequency domain \cite{6} \cite{19}. Short-time Fourier transform (STFT) is a commonly used temporal-frequency analysis method for time series signals, but the time window length is fixed and the resolution of the time domain is relatively low. S-Transform \cite{11} and its derivative version Modified Stockwell Transform (MST)\cite{-1} \cite{5}\cite{24} effectively improve the resolution of the time domain and can effectively extract temporal-frequency features. Discrete Wavelet Transform \cite{20}\cite{26}, Relative Wavelet Energy (RWE) \cite{7}\cite{20} can control the scale by adjusting parameters, can also extract temporal-frequency features, and obtain good results in classification tasks. Those temporal-frequency analysis methods are widely used in ECoG signal processing. In addition, there are many other feature extraction methods, such as Autoregressive Coefficients (AR) \cite{8} Petrotosian Fractal Dimension (PFD) \cite{9} Hjorth Complexity (HC) Hjorth Activity, Mobility, Complexity \cite{10}, common spatial patterns (CSP) \cite{14}\cite{17}\cite{25}, Band Power \cite{18}, conditional random fields (CRF) \cite{15}. The features extracted by these methods are redundant, so the features need to be selected and reduced in size. Commonly used methods include genetic algorithm \cite{13} and principal component analysis (PCA) \cite{20}\cite{26}. The optimized feature subset will be passed to the classifier for final classification. For classifier common methods include support vector machine (SVM), (logistic regression (LR), linear discriminant analysis (LDA) and k-nearest neighbor (kNN), probabilistic neural network (PNN), etc. \\
Traditional methods require manual design the suitable feature extractor, which requires prior knowledge of the relevant parameters. For small data sets, when the similarity between different categories is low and the computational cost is not high, traditional methods can handle it well. However, for large data sets, especially the similarity between data or categories is high, the manually designed features will be not large enough for classification task. For example, the ECoG signal dataset related to visual classification \cite{12} this data has multiple categories, with high similarity cross categories, big data diversity within category, complex visual information. It requires a larger dimensional feature vector to construct an effective feature space which can provide enough discrimination information. At this situation, deep learning methods can show their capability. For temporal series signal processing, LSTM \cite{27} have its reputation, it can effectively extract long-term correlation information. \cite{28} used LSTM to process ECoG signals related to Motion Imaginary and achieved good results in classification task. \cite{29} combined 2D convolution modules with LSTM to prove that the exploited spatial correlation between neighboring electrodes will benefit to predication of the hand trajectory. \cite{12} used vision-based ECoG signals which were collected from the brains of primates. The authors compared the processing capabilities of the convolutional neural network (CNN) and the Bi-LSTM \cite{30} network. The convolutional neural network used two modules, channel-wise transform and channel-level self-attention. The experimental results show that the convolutional neural network effectively extracts the characteristics of the time series signal, which has 6\% higher in the classification accuracy than it of the Bi-LSTM, reaching the SOTA level in the classification task of this vision-based ECoG signal dataset. However, the convolutional neural network and LSTM network in \cite{12} only focus on information in time domain, does not consider the  spatial information of the electrodes. In \cite{0}, which uses the same data set, the authors used MST to perform time-frequency analysis on the ECoG signal to obtain a 3D feature space,  designed a minimalist and explanatory network to extract features from the spatial domain, and used a parallel network to simultaneously process the real and imaginary data. The final performance was significantly higher than \cite{12}. \cite{0} is more like a basic research work,  combining the advantages of traditional methods MST and deep learning methods to explore the characteristics of ECoG signals that can carry visual information. This means that the model inherits the shortcomings of MST, requiring empirical setting of MST parameters, slow temporal-frequency analysis processing, and cannot be integrated into the model network. It requires huge hard disk space to store 3D feature data, this also greatly affects the training speed of the network.\\
Based on the findings of \cite{0}: the frequency domain and spatial domain of ECoG carry more visual information, the effective use of frequency domain and spatial domain features is conducive to improving model performance. In this paper, we propose the Bi-Band ECoGNet model. the Bi-BCWT module in our model has similar functions to MST, mapping the ECoG signal to the 3D feature space, using different length of convolution kernel can effectively take into account both the high-frequency and low-frequency information, improved the feature extraction ability; in addition, we improved the Spatial Filter structure of the model in [0], using a 2D multi-layer structure called spatial-temporal feature encoder, can more effectively extract spatial domain features; finally, the Bi-Band ECoGNet model is more lightweight, has higher performance, and does not require additional space to store data, greatly improving data processing speed. \\
In this paper, introduction of our proposal will be discussed in section II; then model training \& the results comparison in section III; in section IV, we present the ablation experiments; final part is the conclusion.


\begin{figure*}[h!t] 
\centering
\includegraphics[width=7in]{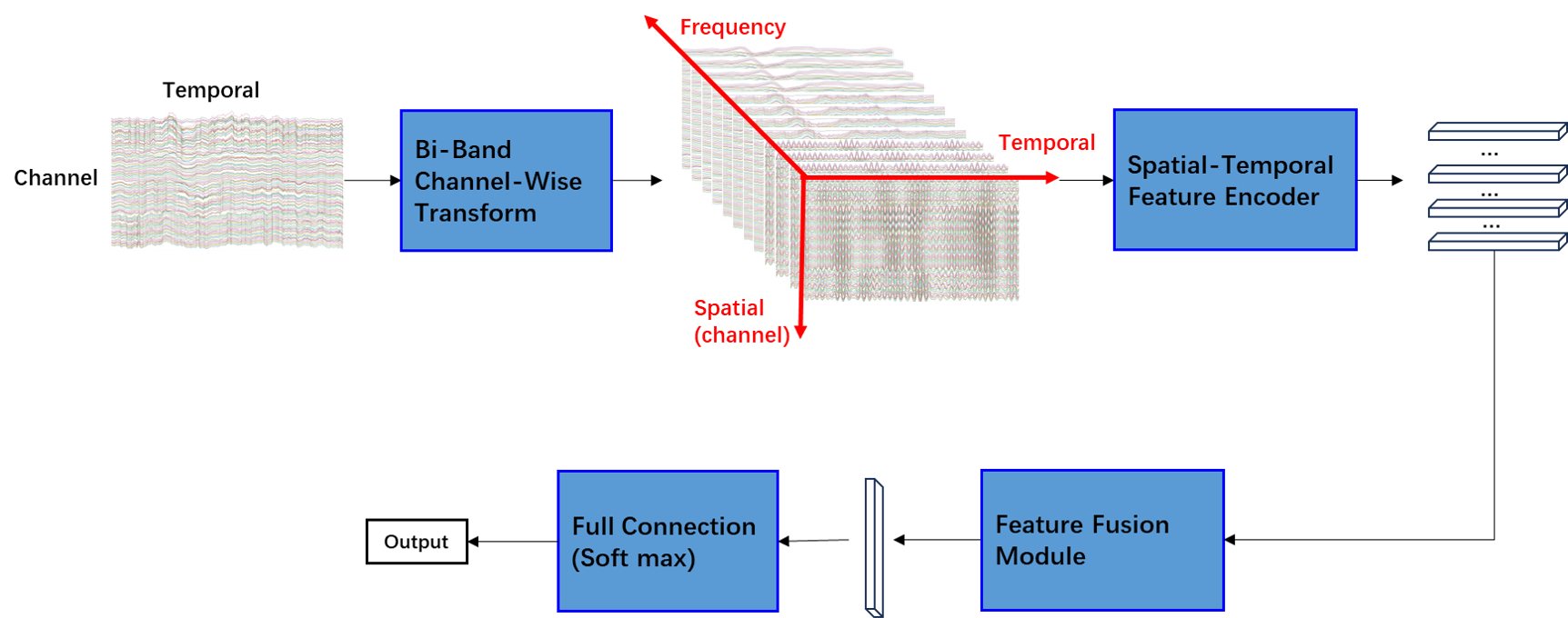}
\caption{Outline of Data Flow of Bi-Band ECoGNet}
\label{fig1:Outline}
\end{figure*}

\section{Method}
Figure \ref{fig1:Outline} shows the entire process of proposed model. There are four main stages, which also correspond to the four parts of the model (blue modules in Figure \ref{fig1:Outline}). The function of Bi-Band channel wise transform is similar to MST \cite{-1}, which extracts the temporal-frequency information of ECoG. After this processing, we can get a 3D feature space. Spatial-Temporal Feature Encoder focuses on the spatial-temporal features in each frequency feature map, extract spatial feature. The following Feature Fusion Module and Full Connection module will do features selection and classification respectively. These two parts of the structure are derived from the structure of EEGNet \cite{32}. The specific structure of each part can be found in Table \ref{table_structure}.

\subsection{Bi-Band Channel-Wise Transform}
Time series signals such as ECoG are non-stationary time series signals, real-time frequency information is one of the important features. There are many methods for temporal-frequency analysis, such as the traditional methods introduced in Section I. However, these methods require manual design of related parameters, which would be tedious and complicated. The convolution layer in the convolutional neural network can provide a solution. We can use the convolution kernel to implement the function of temporal-frequency analysis. In addition, the parameters of the convolution kernel can be obtained through model training, so we can avoid manual design and parameters adjustment.
Here we use the TCN network for time series signals. The structure in TCN uses 1-dimensional convolution layer like \cite{32}. Compared with the TCN module used in \cite{12}, it is smaller and lightweights. The convolution kernel parameters are similar with the parameters of the time window function in the MST method. These parameters can be learned through model training, which greatly reduces the tedious steps of manually adjusting parameters and time consume calculation of MST.\\
Bi-Band Channel-Wise Transform (Bi-BCWT) is composed of 64 TCN small modules, half of which have the kernel length of 512, the other half have kernel length of 32. The two different convolution kernels have different tendencies in the frequency domain. Experimental comparisons show that the 512-convolution kernel pays more attention to low-frequency part, while the 32-convolution kernel is more inclined to relatively high-frequency part. The combination of these two different sizes of convolution kernels can cover a wider range of frequency. The structure of Bi-BCWT is shown in Figure \ref{fig2:bibcwt}:

\begin{figure}[h!t] 
\centering
\includegraphics[width=3.4in]{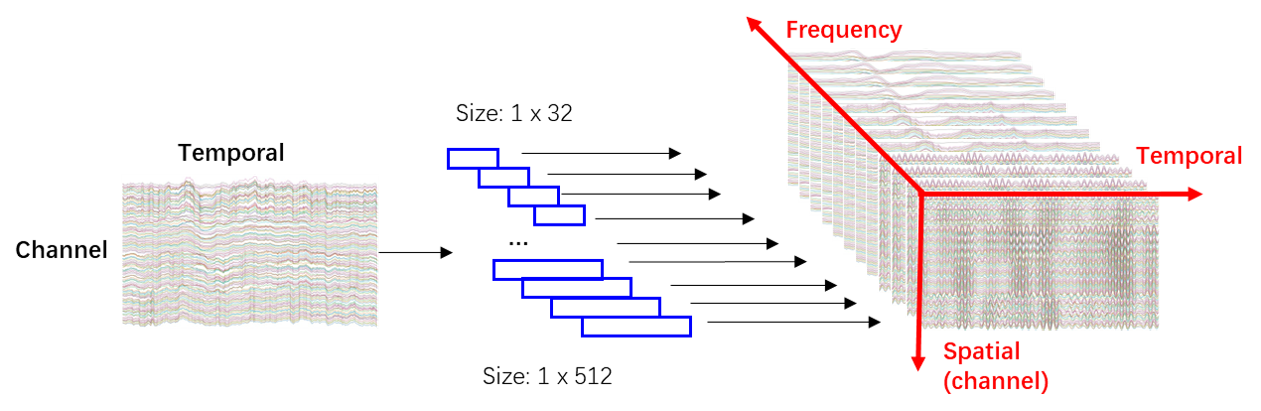}
\caption{\textbf{Structure of Bi-Band Channel-Wise Transform}\\  Each blue rectangle represent one TCN, there are 2 different size of TCNs. short for length of 32, long for length of 512. Each TCN will generate one Spatial-Temporal Feature map}
\label{fig2:bibcwt}
\end{figure}

After the processing of this module, as shown in Figure \ref{fig2:bibcwt}, we can get a 3D feature space (frequency-spatial-temporal). The Bi-BCWT module determines the feature of the frequency domain. In our model, there are 64 Spatial-Temporal feature maps, which correspond to the 64 TCN modules.

\subsection{Spatial-Temporal Feature Encoder}
An 8*16 electrode array was used for ECoG acquisition, so there is corresponding relative spatial information between each electrode. The Spatial Filter structure of MST-ECoGNet \cite{0} uses a 1-dimensional spatial filter, which cannot effectively extract this 2-dimensional spatial feature. Here, in our model, we use a 2-dimensional convolution network. Experiment shows that this 2D structure can more effectively extract spatial features, improve classification accuracy. The detailed structure is shown in Figure \ref{fig3:Spt_Tmp Encoder}:

\begin{figure}[h!t] 
\centering
\includegraphics[width=3.4in]{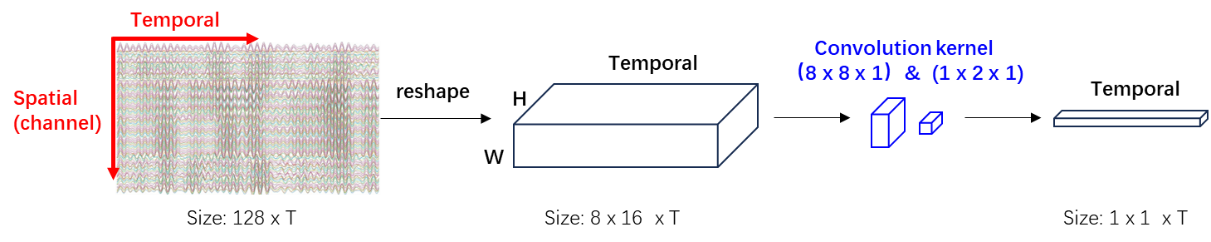}
\caption{\textbf{Structure of Spatial-Temporal Feature Encoder}\\ Each feature map will be reshaped to 3D structure, then processed by 2 layer of 2D kernels, finally output of feature vector. }
\label{fig3:Spt_Tmp Encoder}
\end{figure}

Each spatial-temporal feature map was reshaped into a 3D space (8*16*T), following the same shape as the 8*16 electrode array. Then reshaped feature map will be processed by two layers of convolution, the convolution kernel sizes are 8*8*1 and 1*2*1 respectively. The final output feature vector is in size of 1*1*T. Through this way, the spatial-temporal features are extracted.

\subsection{Feature Fusion Module \& FC}

The fusion module use the same module in EEGNet \cite{32}, detail structure can refer to Table \ref{table_structure}, which is compact and can blend the features of each frequency domain and further reduce the temporal domain. Finally, the FC layer is connected, output a 6-dimensional final classification vector.

\begin{table*}[ht]
\centering
\caption{\textbf{Model Structure Details} N means the number of kernel; G means the 'group' parameter in nn.Cov3d function.}
\label{table_structure}
\resizebox{\textwidth}{!}{
\begin{tabular}{l l l l l l}

\hline
Part &Layer & \multicolumn{2}{l}{[Kernel Size] * N / G} &Output &Reference \\

\hline
\multirow{3}{*}{Bi-BCWT} &nn.Cov3d    &\multicolumn{2}{l}{[1 x 1 x 32] * 32 / 1}  &32 x 128 x 300  &\multirow{2}{*}{Input：ECoG signal of size: 1 x 128 x 300}\\
{} &nn.Cov3d    &\multicolumn{2}{l}{[1 x 1 x 512] * 32 / 1]}  &32 x 128 x 300  &{}\\
{} &Catenate    &\multicolumn{2}{l}{-}  &64 x 128 x 300  &catenate above conv3d output together\\

\hline
\multirow{6}{*}{Spatial-Temporal Encoder} 
    &Reshape    &\multicolumn{2}{l}{-}  &64 x 8 x 16 x 300  &reshape feature map into 3D structure\\
{}  &nn.Cov3d &\multicolumn{2}{l}{[8 x 8 x 1] * 128 / 64} &128 x 1 x 2 x 300    &    \\ 
{}  &nn.Cov3d &\multicolumn{2}{l}{[1 x 2 x 1] * 128 / 128} &128 x 1 x 1 x 300   &    \\ 
{}  &BatchNorm3D &\multicolumn{2}{l}{-} &128 x 1 x 1 x 300   \\
{}  &nn.ELU()    &\multicolumn{2}{l}{-} &128 x 1 x 1 x 300   \\
{}  &nn.AvgPool3d&\multicolumn{2}{l}{[1 x 1 x 4]} &128 x 1 x 1 x 75   \\

\hline{}
\multirow{5}{*}{Feature Fusion}  
    &nn.Cov3d    &\multicolumn{2}{l}{[1 x 1 x 16] * 128 / 128} &128 x 1 x 1 x 75    \\ 
{}  &nn.Cov3d    &\multicolumn{2}{l}{[1 x 1 x 1] * 16 / 1} &16 x 1 x 1 x 75         \\
{}  &BatchNorm3D &\multicolumn{2}{l}{-} &16 x 1 x 1 x 75    \\
{}  &nn.ELU()    &\multicolumn{2}{l}{-} &16 x 1 x 1 x 75    \\
{}  &nn.AvgPool3d&\multicolumn{2}{l}{-} &16 x 1 x 1 x 9     \\

\hline
\multirow{2}{*}{Full Connection}  
    &Flatten   &\multicolumn{2}{l}{-} &1 x 144 &    \\
{}  &nn.lin() &\multicolumn{2}{l}{-} &1 x 6    &    \\

\hline
\end{tabular}
}
\end{table*}


\section{Model Training \& Results}
In the experiment, we used the same visual-based ECoG dataset as \cite{0}\cite{12} to train our proposed model, then compared it with the model proposed in \cite{0}\cite{12} in terms of model size, model performance, and model training speed.

\subsection{Visual-Based ECoG Dataset}
The ECoG data used in experiment was acquired by biology team from Graduate School of Medicine and Dental Sciences, Niigata University. The targets were two macaque monkeys (Subject MonC: 6.1 kg, Subject MonJ: 5.1 kg). The outline of the acquire processing can be seen in Figure \ref{fig4:Outline_record}. For more information, refer to \cite{12} \cite{33}. Here we only give an overview:

\subsubsection{Image set selection}
The images used as visual stimuli were selected from natural images, with a total of 6 categories: building, body part, face, fruit, insect, and tool. After screening by three experts, 1,000 images were finally selected for each category, and the image size was reshaped to 512 x 512;

\subsubsection{ECoG recording}
All animal experiments followed relevant legal requirements, and the experimental subjects were implanted with 8 x 16 ECoG electrode arrays on the surface of the inferior temporal cortex(ITC) of the brain, total 128 electrode channels. In fact, subject MonJ had another implanted electrode array with 64 channels. But here we only utilized 128 channels of data which covered the surface of the ITC. In this way, the data of the two subjects remained consistent.

\subsubsection{Recording Process}
The subjects were trained with a visual fixation task to keep their gazes within ±1.5 degree of the visual angle around the fixation target (diameter: 0.3 degree). Eye movements were captured with an infra-red camera system at a sampling rate of 60 Hz. Stimuli were presented on a 15-inch CRT monitor (NEC, Tokyo, Japan) with a viewing distance of 26 cm. In each trial, after 450 ms of stable fixation, each stimulus was presented for 300 ms, followed by a 600 ms blank interval. The signals were differentially amplified using an amplifier (Plexon, TX, USA or Tucker Davis Technologies, FL, USA) with high- and low-cutoff filters at 300 Hz and 1.0 Hz, respectively. All subdural electrodes were referenced to a titanium screw that was attached directly to the dura at the vertex area. The recording wad was conducted at a sampling rate of 1 kHz per channel.

\begin{figure}[h!t] 
\centering
\includegraphics[width=3.4in]{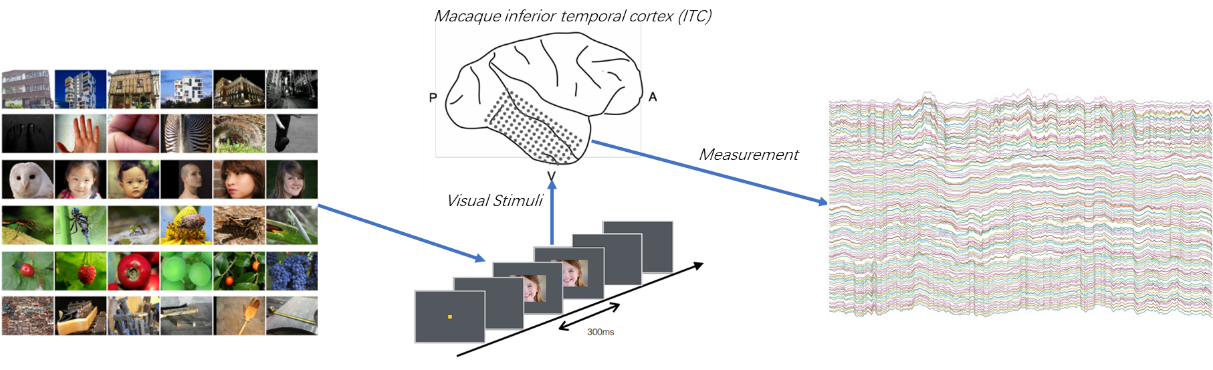}
\caption{\textbf{Outline of ECoG Record} Each image will be used as visual stimuli, last 300 ms, ECoG signal will be measured via electrode array.}
\label{fig4:Outline_record}
\end{figure}

\subsection{ECoG Record Preprocess}
Each ECoG trail contains 3 parts, one is an active interval that lasts for 300 milliseconds. In this interval, picture is displayed to stimulate the monkey's visual system; the other two are the pre-static interval and the post-static interval, picture is not displayed. Please refer to Figure \ref{fig4:Outline_record} in the data acquisition part. For data preprocessing, we select the 300-millisecond interval before the active interval as the background. We calculate the overall mean and variance of the background interval, use them to normalize active interval data. The formula is as follows:

$$E^*_{ac}(n, t)=\frac{E_{ac} (n, t)-m_s}{\sigma_s} $$

$m_s$ $\sigma_s$ indicate overall mean value and standard deviation value of ECoG background respectively. $E_{ac} (n, t) $indicates the ECoG data in active interval, which has $n$ channels, $t$ sampling points. $E^*_{ac} (n, t) $ indicates the ECoG data after processing, this is the input of our model. The purpose of this preprocessing is to exclude the influence of the monkey's state during different experiment period as much as possible. The background interval reflects the basic state of the monkey, if exclude the influence of this period, the remained changes is caused by visual stimulation, which is what we are concerned about.

\subsection{Model Training}

The subject MonC dataset has 21,586 trails, the subject MonJ dataset has 13,445 trails, totaling 35,031 trails of data. Before starting training, the experimental data will be preprocessed. In order to obtain more accurate experimental results, the dataset is evenly divided into 5 parts. In verifying the performance of the model, we use 5 cross fold validation, all result is the average accuracy of 5 times test. The experimental model is completed using pytorch. We trained each model for 400 epochs with a batch size of 128. We used the Adam optimizer \cite{34} for optimizing model parameters, with a learning rate of 1.5e-6, a weight decay of 0.0001, β1 = 0.9, and β2 = 0.999. The operating platform information is as follows: CPU Intel Xeon w5-2455x, memory 64 GB, GPU NVIDIA GeForce RTX 4090 24GB

\subsection{Results \& Comparison}

Table \ref{table_results} shows the comparison results. Compared with the deep learning method \cite{12}, our model is smaller and has better performance. The accuracy of the experimental subject MonC reached 54.15\%, an increase of 17.35\%. The experimental subject MonJ also improved by 8.39\%. In terms of model size, our model is only 0.0396M, which is less than one-tenth of the size of the model \cite{12}. Compared with model in \cite{0}, which combined  traditional method and deep learning method, our model performance is improved by an average of 1.24\% with smaller size. More importantly, training speed of our model is increased by 6 times and no additional storage space is required to save feature data. According to the comparison results, we can see that our model is smaller, higher performance, faster to train and is more suitable for application in the field of BCI.
.

\begin{table*}[ht]
\centering
\caption{\textbf{Model Performance Comparison} compare with previous work \cite{0}\cite{12}, our model has the best performance, faster speed with smaller size}
\label{table_results}
\resizebox{0.75\textwidth}{!}{
\begin{tabular}{|c|c|c|c|c|c|c|c|c|}

\hline
\multicolumn{3}{|c|}{\multirow{2}{*}{Model}} &\multicolumn{2}{c|}{Params / M} &\multicolumn{2}{c|}{Accuracy / \%} &\multicolumn{2}{c|}{Speed / s/epoch}\\
\cline{4-9}

\multicolumn{3}{|l|}{} &{MonJ} &{MonJ} &{MonC} &{MonJ} &{MonC} &{MonJ} \\
\hline

\multicolumn{3}{|l|}{Base Model \cite{12}} &{0.451} &{0.451} &{36.8±0.51} &{27.59±0.73} &{-} &{-} \\
\hline
\multicolumn{3}{|l|}{MST-ECoGNet \cite{0}}  &{0.0396} &{0.0488} &{53.43±0.55} &{34.22±0.78} &{189.0} &{139.0}\\
\hline

\multicolumn{3}{|c|}{Bi-Band ECoGNet}  &\textbf{0.0396} &\textbf{0.0396} &{\textbf{54.15±1.15}} &\textbf{35.98±1.15} &{\textbf{33.0}} &{\textbf{21.2} }\\

\hline

\end{tabular}
}
\end{table*}

\section{Research Experiments}
The ECoG data in this experiment is unique. In addition, unlike the Motion Imagery dataset, this ECoG dataset is collected from the Inferior Temporal Cortex. The Inferior Temporal (IT) Cortex is the cerebral cortex on the inferior convexity of the temporal lobe in primates. It is crucial for visual subject recognition and considered as the final stage in the ventral cortical visual system \cite{34}. Few studies have been conducted on the classification task based on vision-based ECoG. Therefore, on the basis of the proposed model, we conducted ablation experiments on the key modules and parameters of the model to ensure that we understood the functions of the modules and improved the reliability of the model. In addition, in the channel importance test, we have new findings that are worth sharing with experts in the biological team. In the ablation experiment, we did not distinguish the data set of different experimental subjects and mixed them together for testing.

\subsection{Number of TCNs \& Length of Kernel}

\begin{figure}[h!t] 
\centering
\includegraphics[width=3.4in]{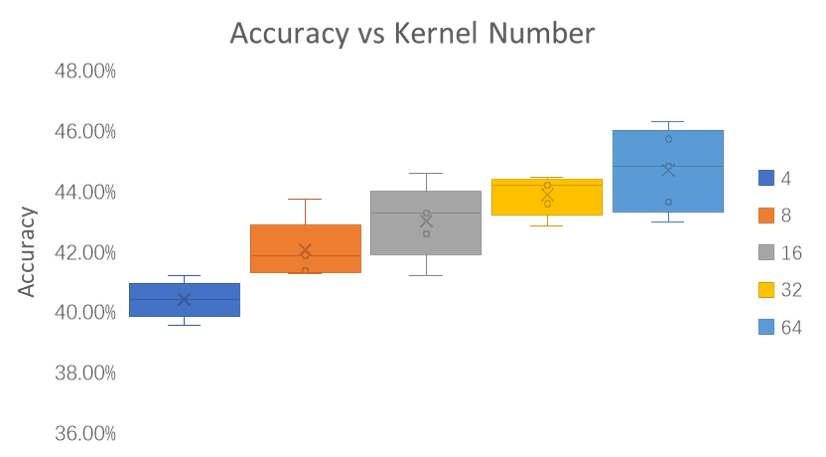}
\caption{\textbf{Performance of Different Kernel Number} Results show that the more kernel (TCN), the better performance.}
\label{fig5:num_tcn}
\end{figure}

\begin{figure}[h!t] 
\centering
\includegraphics[width=3.4in]{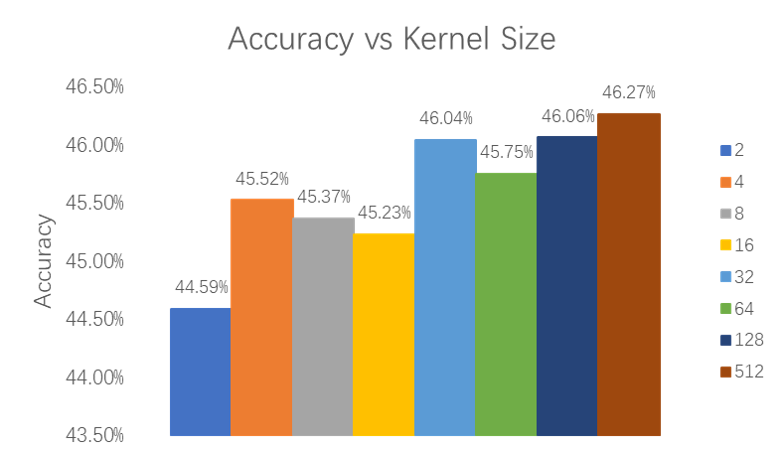}
\caption{\textbf{Performance of Different Size of Kernel} The size (length) of kernel has influence on performance. Finally 32 \& 512 kernel size are selected, base on the performance.}
\label{fig6:len_tcn}
\end{figure}

The important component of Bi-BCWT is the TCN module. TCNs number and the length of its kernel will determine the performance of Bi-BCWT. The number of TCNs determines the number of frequency feature maps; the length of kernel affects the range of the frequency domain. We test the model performance using different number of TCNs and different length of the convolution kernel. The experimental results are shown in Figures \ref{fig5:num_tcn} and Figure \ref{fig6:len_tcn}. From comparison results, we can see that: 1) as the number of TCNs increases, the classification accuracy increases; 2) Kernel Length also affects the performance of the model, it shows relatively high performance when the kernel length is 512 and 32. \\Since the main function of TCN is to extract frequency information of each channel, frequency is determined by the kernel parameters. The increase in the number of TCNs also increases the diversity of kernel frequencies. In other words, it will determine the frequency resolution, the larger the number, the higher the frequency resolution. The size of the kernel length will affect the model's choice of high and low frequency. A larger length is more inclined to the low frequency band. On the contrary, a smaller length will have a greater chance of learning relatively high frequency domain.

\begin{figure}[h!t] 
\centering
\includegraphics[width=3.4in]{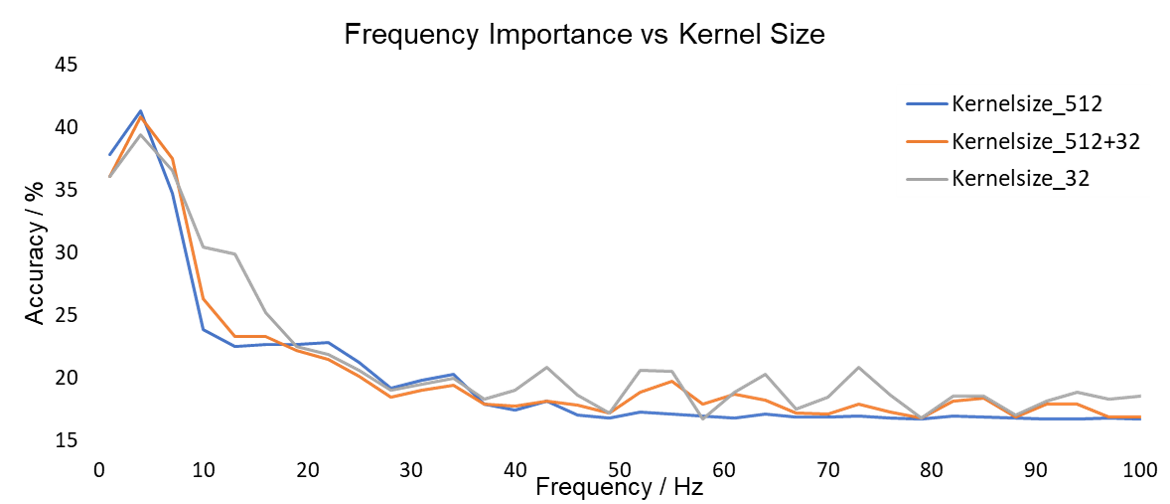}
\caption{\textbf{Frequency Importance of Different Kernel Size} Results show that model with big size of kernel such as 512, main focus on low frequency information, small size kernel like 32 can get information from relatively high frequency. However, when combine these two kernel size together, both low frequency and high frequency information can provide contribution, so this combination can use both low and high frequency information.}
\label{fig7:frq_len_tcn}
\end{figure}

To prove this conjecture, we designed a verification experiment: designing an FIR band-pass digital filter with bandwidth of 5Hz, the cutoff frequency gradually increases from 1Hz to 300Hz, with a step size of 3Hz. For each step, original ECoG filtered by the filter, then fed into pretrained model for accuracy test. Detail algorithm refer to Algorithm \ref{alg:alg1_frq}. The experimental results are shown in Figure \ref{fig7:frq_len_tcn}: The model with length 512 shows high performance in the low frequency range. On the contrary, the model with length 32 can extract relatively high frequency information.

\begin{algorithm}[h!t]
    \caption{Frequency Importance Test}
    \label{alg:alg1_frq}
    \KwIn{original ECoG dataset $E_{ori}$}
    \KwOut{ model accuracy list $Acc=[acc_0,\ldots,acc_i,\ldots,acc_n]$,  $acc_i$ means model accuracy after $i^{th}$ band pass filter process original ECoG dataset }
    \# Each band pass filter has bandwidth of 5Hz\\
    initialization\\
    $Acc$=$[ ]$   \qquad \qquad \qquad \qquad \# accuracy list \\ 
    $fr_{on}$=$1Hz$  \qquad \qquad \,  \qquad \# start frequency\\
    $fr_{off}$= $fr_{on}$+$5Hz$   \quad \qquad\# cutoff frequency  
    
    \While{$fr_{on}< $ 300Hz}{
        \# Implement of filter function via package signal\;

        $filt_i=signal.firls(fr_{on} ,fr_{off}) $\\
	$E^*_i=filt_i (E_{ori})	$		\\
	$acc_i=Model(E^*_i)$   			\\
	$Acc.append(acc_i)$				\\
	$fr_{on}=fr_{on}+3Hz$			

       }
\end{algorithm}

\subsection{Bi-Band TCNs vs Solo-Band TCNs}
Based on the findings in \textbf{Experiment A}, different kernel length leads different performances in different frequency domain. Larger kernel length tends to low frequency, smaller length trend to relative higher frequency. In order to extract both low-frequency and relative high-frequency information, mixed TCNs of different lengths was proposed. Keeping the number of TCNs unchanged, half used length of 32, half used length of 512. Comparison experiments refer to Table \ref{table_ablation}, verified that the performance of the model was improved after using mixed TCNs. Frequency test also showed that the model did can extract larger range frequency information, both high and low frequency, refer to Figure \ref{fig7:frq_len_tcn}.

\subsection{Electrodes Array Spatial Domain}
\begin{figure}[h!t] 
\centering
\includegraphics[width=3.2in]{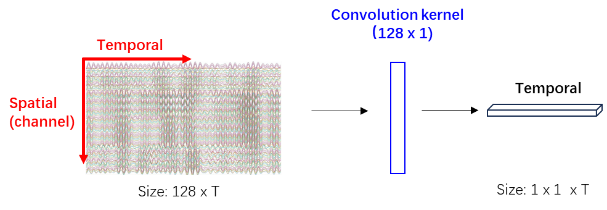}
\caption{\textbf{Spatial-Temporal Feature Encoder of 2DECoGNet} Each Spatial-Temporal feature map will processed by 1D convolution kernel, then get feature vector, the big difference compare with 3DECoGNet is that no reshape operation on the feature map.}
\label{fig9:2decognet}
\end{figure}

Work \cite{0} mentioned that the spatial domain features of ECoG are helpful for classification tasks, 1D convolution kernel was used to extract spatial domain features. However, ECoG acquisition used an 8 x 16 rectangular electrode array. This made the channel of ECoG have 2D spatial domain. To better extract spatial features, Spatial-Temporal feature Encoder was proposed. As introduced previously, the spatial-temporal feature map was reshaped into 3D data (8 x 16 x T), then two layers of 2D convolution kernels were used to extract features. To verify its performance, we conducted a comparison experiment:\\
1)	2DECoGNet: same structure used in [0], keep spatial-temporal feature map as 2D data structure, directly use 1D convolution kernel extract spatial features, the structure is shown in Figure \ref{fig9:2decognet};\\
2)	3DECoGNet: the network with Spatial-Temporal feature Encoder.

The comparison results are shown in Table \ref{table_ablation}. Compared with 2DECoGNet, 3DECoGNet has higher performance with accuracy improved by 1.04\%. The special electrode array used in ECoG acquisition contains 2D spatial information that can help improve the accuracy of multi-classification tasks.

\begin{table}[ht]
\centering
\caption{\textbf{Ablation Study on Kernels length \& Spatial-Temporal Encoder} All test on mix dataset of 2 subjects. Accuracy uses 5 cross validation average value.}
\label{table_ablation}
\resizebox{0.48\textwidth}{!}{
\begin{tabular}{c c c c c c c c}

\hline
\multicolumn{4}{c|}{Kernel Length} &\multicolumn{2}{c}{Spatial-Temporal Encoder} &\multirow{2}{*}{Model Size / M} &\multirow{2}{*}{Accuracy / \%}\\
\cline{1-6}

256 &512    &32  &512+32 &2D Encoder    &3D Encoder\\
\hline

\ding{52}  &   &   &   &{\ding{52}} &     &0.0380    &45.09   \\
\ding{52}  &   &   &   &   &{\ding{52}}   &0.0305    &46.13   \\
\hline

{}  &\ding{52}   &   &   &   &{\ding{52}}   &0.0648    &46.27   \\
{}  &   &\ding{52}   &   &   &{\ding{52}}   &0.0161    &46.04   \\
{}  &   &   &\ding{52}   &   &{\ding{52}}   &0.0396    &46.77   \\

\hline
\end{tabular}
}
\end{table}

\subsection{Individual Channel Influence}
\begin{figure}[h!t] 
\centering
\includegraphics[width=3.2in]{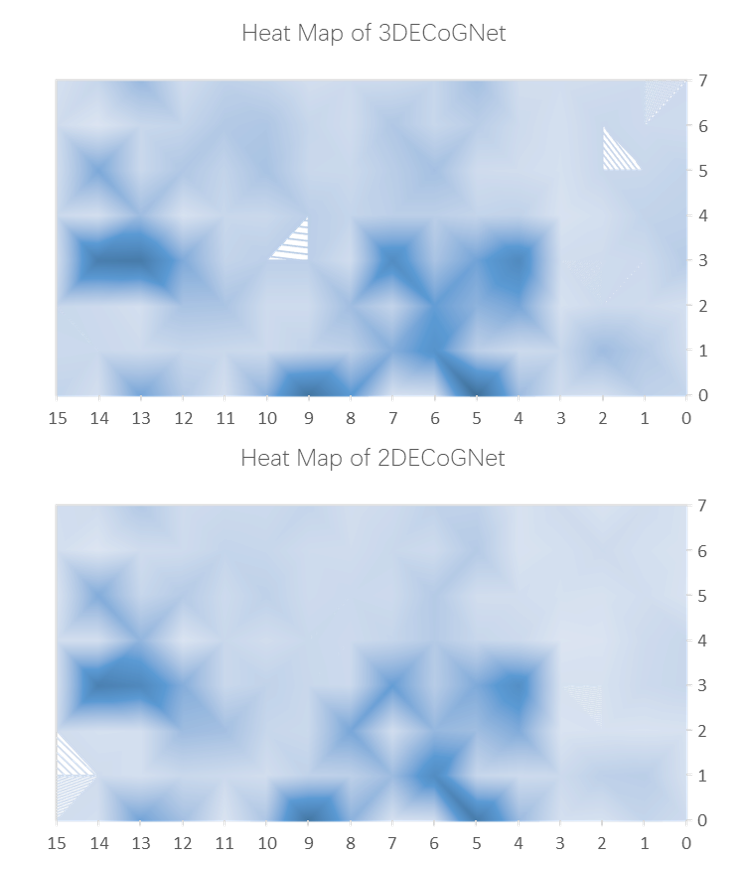}
\caption{\textbf{Heat Map for Each Channel of 3DECoGNet \& 2DECoGNet} 128 channel reshaped to 8 x 16, which is same as electrode array, number of x y axis is used to define the position of electrode. Deeper color means bigger value which indicates more important for model performance. Results show that 3DECoGNet can extract more information from more channels, such as (10, 5), (6, 5), (2, 1)etc. compared with 2DECoGNet heat map, }
\label{fig11:3d_heatmap}
\end{figure}

The above \textbf{Experiment C} shows that after adopting the 3D data structure, we can clearly see the improvement in accuracy. To find out the reason for the improvement, we focus on the signal channel and wonder if 3DECoGNet can learn channel information better. To this end, we designed a comparison experiment to see the contribution level of each channel. The experimental method is very simple:\\
ECoG data composes of 128 channels, each time only one channel data are set to 0.  The modified ECoG dataset is fed into the pre-trained model to obtain the accuracy. This accuracy will decrease. Compare with the accuracy of original ECoG dataset, we can get the value of decrease. This value can be used as contribution indicator of the channel which is set to 0. Repeat this procedure until all 128 channels have been test. Detail algorithm refer to Algorithm \ref{alg:alg2_chn}. Finally, all channel contribution level are drawn on the 8 x 16 plane map, then we can get the heat map.\\
Compare the heat maps of 3DECoGNet with  2DECoGNet Figure \ref{fig11:3d_heatmap}. We can find that the heat region of 3DECoGNet is larger than that of 2DECoGNet, for example, electrodes (10, 5), (6, 5) and (2, 1) etc. Therefore, we can conclude that the 2D kernel helps to extract more spatial information, which is conducive to the multi-classification task of ECoG. Please note that for this test we mixed the two monkeys' dataset, this is little different with \textbf{Experiment E}.

\begin{algorithm}[h!t]
    \caption{Channel Test}
    \label{alg:alg2_chn}
    \KwIn{original ECoG dataset $E_{ori}$}
    \KwOut{accuracy change list $D=[d_0,\ldots,d_i,\ldots,d_n]$,  $d_i$ means model accuracy change after $i^{th}$ channel value set to 0}
    \# total channel number 128\\
    initialization\\
    $D$=$[ ]$   \qquad \qquad \qquad \qquad \qquad \, \# accuracy change list \\ 
    $Acc_{ori}$=$Model(E_{ori})$  \qquad \qquad  \# original accuracy\\
    $i$= $0$  
    
    \While{$i < $ 128}{
        
        Set $i^{th}$  channel of $E_{ori}$ to 0, get modified ECoG dataset $E^*_i$\\
	$acc_i=Model(E^*_i)$   			\\
        $d_i=Acc_{ori} - acc_i$         \\	
	$D.append(d_i)$				\\
	$i+=1$			
       }
    \# final delta list reshapes to 8 x 16 shape
\end{algorithm}

\subsection{Difference between Subject}

\begin{table}[ht]
\centering
\caption{\textbf{Subject-Cross Performance Comparison} Results show that model that trained and test on different monkey dataset, show big difference. when we train model on one monkey dataset, test on other monkey dataset, it show very bad performance at random guess level.}

\label{table_crosssubject}
\resizebox{0.3\textwidth}{!}{
\begin{tabular}{|c|c|c|c|}

\hline
\multicolumn{2}{|c|}{\multirow{2}{*}{Dataset/Accuracy}} &\multicolumn{2}{c|}{Test} \\
\cline{3-4}
\multicolumn{2}{|l|}{}  &MonC   &MonJ \\
\hline
\multirow{2}{*}{Train}  &MonC   &54.15\%  &16.76\%  \\ 
\cline{2-4}
{}                      &MonJ   &15.88\%  &35.98\%  \\ 
\hline
\end{tabular}
}
\end{table}

\begin{figure}[h!t] 
\centering
\includegraphics[width=3.2in]{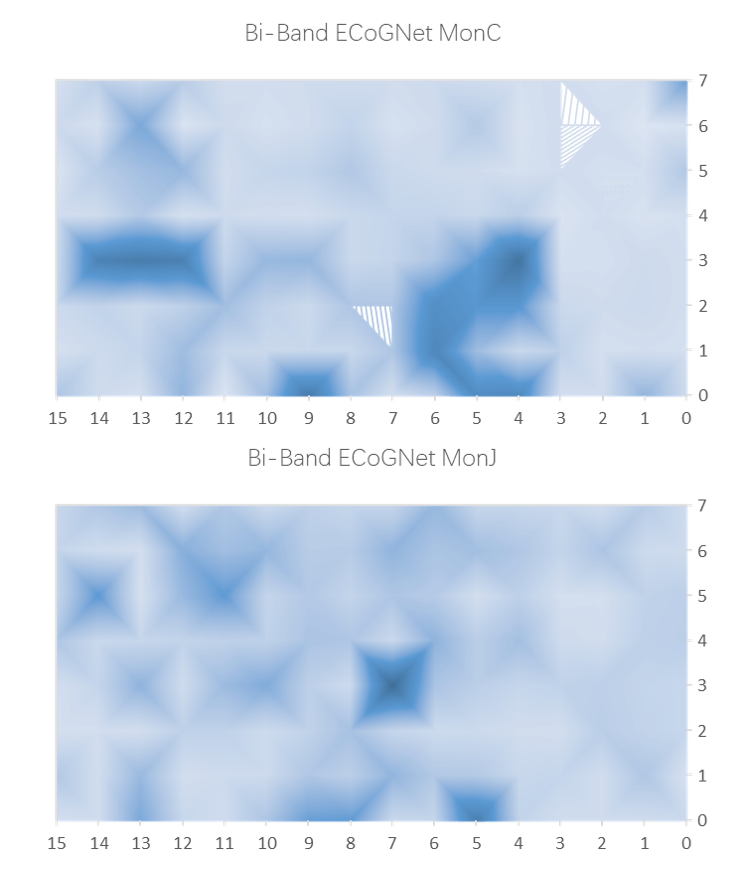}
\caption{\textbf{Heat Map of MonC \& MonJ Dataset} Results show that there is big difference in channel heat map between MonC \& MonJ dataset.The position where visual information comes from is quite different between 2 subjects.}
\label{fig15:chn_monC}
\end{figure}

From \cite{0}\cite{12} we can see the big difference between the model performance when the model was trained on each subject dataset separately. The same difference is also found in proposed model, with a difference of 19.17\%. A subject-cross test was conducted: model trained on one subject`s dataset then test on the other subject`s dataset. \\
Table \ref{table_crosssubject} shows the cross test results: model trained on one subject dataset can not get good accuracy on other subject dataset. This means each dataset has its own specific. The possible reason for this specific would be biology difference between these 2 monkeys. \\
Following this view, try to verify the biology difference, channel importance test was conducted on each subject dataset. Used same method in \textbf{Experiment C} refer to algorithm \ref{alg:alg2_chn}. The heat maps of the two subjects are shown in Figure \ref{fig15:chn_monC}. Compare with two heat map, there is a big difference. Different channel heat map means visual signal comes from different position. Considering same electrode array, same position, similar subject, same visual task, these difference would indicate that special pattern exists in visual processing of brain, different cluster of cortex cell activate during that processing. Different individual process visual information in its special pattern. However, above hypothesis is based on the data analysis results, we are not experts in biology, we would like to share these findings and discuss with the biology experts in Niigata University.

\section{Conclusion}
There are very few studies on ECoG signals generated by the visual system. As far as we know, we are the first team to have such dataset. These experimental data are acquired by the excellent biological team at Niigata University. Their efforts have given us the opportunity to peek into the function of the brain's visual system. Although there is still a long way to go to fully understand the full picture of the brain's function, this valuable dataset indeed opened the door for us. Based on the previous work \cite{0}, we further analyzed this dataset and proposed Bi-Band ECoGNet for the multi-classification task of visual ECoG. The Bi-BCWT module in proposed model covers both high-frequency and low-frequency features, can effectively map ECoG signals to 3D feature space, greatly improving the efficiency of data processing compared to the MST method in \cite{0}. In addition, the Spatial-Temporal Feature Encoder uses two layers of 2D convolution kernels, which is more conducive to extracting spatial domain features than the 1-dimensional Spatial Filter in \cite{0}. Overall, the proposed model has better performance than MST-ECoGNet, accuracy improved by 1.24\%, and what`s more important, our model training speed is 6 times faster, no need additional resource for dataset.\\
Furthermore, in the ablation experiment, we systematically analyzed how the model's hyper parameters affect the feature space, and how the feature impact the accuracy of classification task.\\
1) Frequency importance test: For ECoG multi-classification tasks, distinguish information comes from the low-frequency region, while there is almost no effective information in the high-frequency domain (the accuracy of the high-frequency domain is close to the level of random guessing). The convolution kernel size of Bi-BCWT's TCN determines the preference for the frequency information. Larger convolution kernel tends to low-frequency information, while smaller convolution kernel focus on relatively high-frequency information. The combination of different size of convolution kernel enable the network to take into account both low-frequency and high-frequency information, thereby improving performance.\\
2) The electrode array for collecting ECoG data is a regular rectangle, which makes the ECoG channels also have unique 2D spatial information. 2D convolution kernel in Spatial-Temporal Feature Encoder can more effectively extract spatial information, so that make our model get a 0.5\% performance improvement.\\
3) The contribution of each ECoG channel to the task is different. In the experiment, we found that only small part of channels can provide effective information, showing a relatively concentrated characteristic. Effective selection of channels will be the next research topic, as mentioned in \cite{35}\cite{36}.
4) Analyzing the difference between subjects, compare the channel heat map between 2 subject, we can find there is a big difference in the heat map. Base on this finding, we guess that there exist some special biology difference in each subject. Even the same visual task, different individual process that visual information in different pattern. This hypothesis is worth to share with biological experts in Niigata University.


\section*{Acknowledgment}
The author of this paper is financially support by CSC (China Scholarship Committee) and RA of OKATANI Lab in GSIS of Tohoku University.


\bibliographystyle{IEEEtran}
\bibliography{IEEEabrv,biblio}

\begin{thebibliography}{10}
\providecommand{\url}[1]{#1}
\csname url@samestyle\endcsname
\providecommand{\newblock}{\relax}
\providecommand{\bibinfo}[2]{#2}
\providecommand{\BIBentrySTDinterwordspacing}{\spaceskip=0pt\relax}
\providecommand{\BIBentryALTinterwordstretchfactor}{4}
\providecommand{\BIBentryALTinterwordspacing}{\spaceskip=\fontdimen2\font plus
\BIBentryALTinterwordstretchfactor\fontdimen3\font minus \fontdimen4\font\relax}
\providecommand{\BIBforeignlanguage}[2]{{%
\expandafter\ifx\csname l@#1\endcsname\relax
\typeout{** WARNING: IEEEtran.bst: No hyphenation pattern has been}%
\typeout{** loaded for the language `#1'. Using the pattern for}%
\typeout{** the default language instead.}%
\else
\language=\csname l@#1\endcsname
\fi
#2}}
\providecommand{\BIBdecl}{\relax}
\BIBdecl

\bibitem{0}
\BIBentryALTinterwordspacing
JI,~C., ``Explainable mst-ecognet decode visual information from ecog signal,'' 2024. [Online]. Available: \url{https://arxiv.org/abs/2411.16165}
\BIBentrySTDinterwordspacing

\bibitem{-1}
Jiao,~B., Rui,~Y., Gao,~M., Fei,~H., and Yu,~Q., ``A modified s transform with adjustable window function,'' in \emph{2019 IEEE International Conference on Signal, Information and Data Processing (ICSIDP)}, 2019, pp. 1--5.

\bibitem{1}
\BIBentryALTinterwordspacing
Romanelli,~P., Piangerelli,~M., Ratel,~D., Gaude,~C., Costecalde,~T., Puttilli,~C., Picciafuoco,~M., Benabid,~A., and Torres,~N., ``A novel neural prosthesis providing long-term electrocorticography recording and cortical stimulation for epilepsy and brain-computer interface,'' \emph{Journal of Neurosurgery}, vol. 130, no.~4, pp. 1166 -- 1179, 2019. [Online]. Available: \url{https://thejns.org/view/journals/j-neurosurg/130/4/article-p1166.xml}
\BIBentrySTDinterwordspacing

\bibitem{2}
\BIBentryALTinterwordspacing
Wakuya,~M., Inoue,~T., Imoto,~H., Maruta,~Y., Nomura,~S., Suzuki,~M., and Yamakawa,~T., ``Epileptic seizure–related changes in electrocorticogram, cortical temperature, and cerebral hemodynamics obtained via an implantable multimodal multichannel probe during preoperative monitoring: illustrative case,'' \emph{Journal of Neurosurgery: Case Lessons}, vol.~3, no.~10, p. CASE21694, 2022. [Online]. Available: \url{https://thejns.org/caselessons/view/journals/j-neurosurg-case-lessons/3/10/article-CASE21694.xml}
\BIBentrySTDinterwordspacing

\bibitem{3}
\BIBentryALTinterwordspacing
Islam,~M.~R., Zhao,~X., Miao,~Y., Sugano,~H., and Tanaka,~T., ``Epileptic seizure focus detection from interictal electroencephalogram: a survey,'' \emph{Cognitive Neurodynamics}, vol.~17, p.~1, 2 2023. [Online]. Available: \url{https://www.ncbi.nlm.nih.gov/pmc/articles/PMC9871145/}
\BIBentrySTDinterwordspacing

\bibitem{31}
\BIBentryALTinterwordspacing
Yao,~L., Baker,~J.~L., Schiff,~N.~D., Purpura,~K.~P., and Shoaran,~M., ``Predicting task performance from biomarkers of mental fatigue in global brain activity,'' \emph{Journal of Neural Engineering}, vol.~18, no.~3, p. 036001, mar 2021. [Online]. Available: \url{https://dx.doi.org/10.1088/1741-2552/abc529}
\BIBentrySTDinterwordspacing

\bibitem{4}
\BIBentryALTinterwordspacing
Keene,~D.~L., Roberts,~D., Splinter,~W.~M., Higgins,~M., and Ventureyra,~E., ``Alfentanil mediated activation of epileptiform activity in the electrocorticogram during resection of epileptogenic foci,'' \emph{Can. J. Neurol. Sci}, vol.~24, pp. 37--39, 1997. [Online]. Available: \url{https://doi.org/10.1017/S0317167100021065}
\BIBentrySTDinterwordspacing

\bibitem{6}
Zhou,~F. X. W. Z. D. S. Q. Y.~W., ``Decoding spectro-temporal representation for motor imagery recognition using ecog-based brain-computer interfaces,'' \emph{JIN}, vol.~19, no.~2, pp. 259--272, 2020.

\bibitem{19}
\BIBentryALTinterwordspacing
Liang,~N. and Bougrain,~L., ``Decoding finger flexion from band-specific ecog signals in humans,'' \emph{Frontiers in Neuroscience}, vol.~6, 2012. [Online]. Available: \url{https://www.frontiersin.org/journals/neuroscience/articles/10.3389/fnins.2012.00091}
\BIBentrySTDinterwordspacing

\bibitem{11}
Thangaraj,~K., Muruganandham,~J., Selvaumar,~S., and Jagan,~R., ``Analysis of harmonics using s-transform,'' in \emph{2016 International Conference on Emerging Trends in Engineering, Technology and Science (ICETETS)}, 2016, pp. 1--5.

\bibitem{5}
Jiao,~B., Rui,~Y., Gao,~M., Fei,~H., and Yu,~Q., ``A modified s transform with adjustable window function,'' in \emph{2019 IEEE International Conference on Signal, Information and Data Processing (ICSIDP)}, 2019, pp. 1--5.

\bibitem{24}
\BIBentryALTinterwordspacing
Assous,~S. and Boashash,~B., ``Evaluation of the modified s-transform for time-frequency synchrony analysis and source localisation,'' \emph{eurasipjournals}, 2012. [Online]. Available: \url{http://asp.eurasipjournals.com/content/2012/1/49}
\BIBentrySTDinterwordspacing

\bibitem{20}
Zhao,~H.-b., Yu,~C.-y., Liu,~C., and Wang,~H., ``Ecog-based brain-computer interface using relative wavelet energy and probabilistic neural network,'' in \emph{2010 3rd International Conference on Biomedical Engineering and Informatics}, vol.~2, 2010, pp. 873--877.

\bibitem{26}
\BIBentryALTinterwordspacing
Shi,~H., Yu,~P., and Li,~H., ``The finger flexion related feature extraction method based on wavelet time-frequency analysis in ecog signals,'' in \emph{Proceedings of the 5th International Conference on Computer Science and Application Engineering}, ser. CSAE '21.\hskip 1em plus 0.5em minus 0.4em\relax New York, NY, USA: Association for Computing Machinery, 2021. [Online]. Available: \url{https://doi.org/10.1145/3487075.3487099}
\BIBentrySTDinterwordspacing

\bibitem{7}
Zhao,~H.-b., Yu,~C.-y., Liu,~C., and Wang,~H., ``Ecog-based brain-computer interface using relative wavelet energy and probabilistic neural network,'' in \emph{2010 3rd International Conference on Biomedical Engineering and Informatics}, vol.~2, 2010, pp. 873--877.

\bibitem{8}
Hammon,~P.~S. and de~Sa,~V.~R., ``Preprocessing and meta-classification for brain-computer interfaces,'' \emph{IEEE Transactions on Biomedical Engineering}, vol.~54, no.~3, pp. 518--525, 2007.

\bibitem{9}
\BIBentryALTinterwordspacing
Xu,~F., Zhou,~W., Zhen,~Y., Yuan,~Q., and Wu,~Q., ``Using fractal and local binary pattern features for classification of ecog motor imagery tasks obtained from the right brain hemisphere,'' \emph{International Journal of Neural Systems}, vol.~26, no.~06, p. 1650022, 2016, pMID: 27255798. [Online]. Available: \url{https://doi.org/10.1142/S0129065716500222}
\BIBentrySTDinterwordspacing

\bibitem{10}
Yao,~L. and Shoaran,~M., ``Enhanced classification of individual finger movements with ecog,'' in \emph{2019 53rd Asilomar Conference on Signals, Systems, and Computers}, 2019, pp. 2063--2066.

\bibitem{14}
Jiang,~T., Jiang,~T., Wang,~T., Mei,~S., Liu,~Q., Li,~Y., Wang,~X., Prabhu,~S., Sha,~Z., and Ince,~N.~F., ``Characterization and decoding the spatial patterns of hand extension/flexion using high-density ecog,'' \emph{IEEE Transactions on Neural Systems and Rehabilitation Engineering}, vol.~25, no.~4, pp. 370--379, 2017.

\bibitem{17}
Li,~Y., Koike,~Y., and Sugiyama,~M., ``A framework of adaptive brain computer interfaces,'' in \emph{2009 2nd International Conference on Biomedical Engineering and Informatics}, 2009, pp. 1--5.

\bibitem{25}
Deng,~X., Li,~D., Mi,~J., Gao,~F., Chen,~Q., Wang,~J., and Liu,~R., ``Motor imagery ecog signal classification using sparse representation with elastic net constraint,'' in \emph{2018 IEEE 7th Data Driven Control and Learning Systems Conference (DDCLS)}, 2018, pp. 44--49.

\bibitem{18}
Chong,~L., Hai-bin,~Z., Chun-sheng,~L., and Hong,~W., ``Classification of ecog signals for motor imagery tasks,'' in \emph{2010 2nd International Conference on Signal Processing Systems}, vol.~3, 2010, pp. V3--185--V3--188.

\bibitem{15}
\BIBentryALTinterwordspacing
Saa,~J. F.~D., de~Pesters,~A., and Cetin,~M., ``Asynchronous decoding of finger movements from ecog signals using long-range dependencies conditional random fields,'' \emph{Journal of Neural Engineering}, vol.~13, no.~3, p. 036017, may 2016. [Online]. Available: \url{https://dx.doi.org/10.1088/1741-2560/13/3/036017}
\BIBentrySTDinterwordspacing

\bibitem{13}
\BIBentryALTinterwordspacing
Jain,~R., Jaiman,~P., and Baths,~V., ``Feature engineering for an efficient motor related ecog bci system,'' \emph{bioRxiv}, 2023. [Online]. Available: \url{https://www.biorxiv.org/content/early/2023/04/15/2023.04.01.535201}
\BIBentrySTDinterwordspacing

\bibitem{12}
Date,~H., Kawasaki,~K., Hasegawa,~I., and Okatani,~T., ``Deep learning for channel-agnostic brain decoding across multiple subjects,'' in \emph{2020 8th International Winter Conference on Brain-Computer Interface (BCI)}, 2020, pp. 1--6.

\bibitem{27}
\BIBentryALTinterwordspacing
Zia,~T. and Zahid,~U., ``Long short-term memory recurrent neural network architectures for urdu acoustic modeling,'' \emph{Int. J. Speech Technol.}, vol.~22, no.~1, p. 21–30, mar 2019. [Online]. Available: \url{https://doi.org/10.1007/s10772-018-09573-7}
\BIBentrySTDinterwordspacing

\bibitem{28}
\BIBentryALTinterwordspacing
Du,~A., Yang,~S., Liu,~W., and Huang,~H., ``Decoding ecog signal with deep learning model based on lstm,'' \emph{TENCON 2018 - 2018 IEEE Region 10 Conference}, pp. 0430--0435, 2018. [Online]. Available: \url{https://api.semanticscholar.org/CorpusID:67874163}
\BIBentrySTDinterwordspacing

\bibitem{29}
\BIBentryALTinterwordspacing
Śliwowski,~M., Martin,~M., Souloumiac,~A., Blanchart,~P., and Aksenova,~T., ``Decoding ecog signal into 3d hand translation using deep learning,'' \emph{Journal of Neural Engineering}, vol.~19, no.~2, p. 026023, mar 2022. [Online]. Available: \url{https://dx.doi.org/10.1088/1741-2552/ac5d69}
\BIBentrySTDinterwordspacing

\bibitem{30}
Schuster,~M. and Paliwal,~K., ``Bidirectional recurrent neural networks,'' \emph{IEEE Transactions on Signal Processing}, vol.~45, no.~11, pp. 2673--2681, 1997.

\bibitem{32}
\BIBentryALTinterwordspacing
Lawhern,~V.~J., Solon,~A.~J., Waytowich,~N.~R., Gordon,~S.~M., Hung,~C.~P., and Lance,~B.~J., ``Eegnet: a compact convolutional neural network for eeg-based brain–computer interfaces,'' \emph{Journal of Neural Engineering}, vol.~15, no.~5, p. 056013, jul 2018. [Online]. Available: \url{https://dx.doi.org/10.1088/1741-2552/aace8c}
\BIBentrySTDinterwordspacing

\bibitem{33}
Date,~H., Kawasaki,~K., Hasegawa,~I., and Okatani,~T., ``Deep learning for natural image reconstruction from electrocorticography signals,'' in \emph{2019 IEEE International Conference on Bioinformatics and Biomedicine (BIBM)}, 2019, pp. 2331--2336.

\bibitem{34}
Gross,~C.~G., ``Inferior temporal cortex,'' \emph{Scholarpedia}, vol.~3, p. 7294, 12 2008.

\bibitem{35}
Wei,~Q. and Tu,~W., ``Channel selection by genetic algorithms for classifying single-trial ecog during motor imagery,'' in \emph{2008 30th Annual International Conference of the IEEE Engineering in Medicine and Biology Society}, 2008, pp. 624--627.

\bibitem{36}
Y,~L. X. W. D. Z. B. F. C. C. J. X. M.~C., ``A review on electroencephalogram based channel selection,'' \emph{Sheng Wu Yi Xue Gong Cheng Xue Za Zhi}, vol.~41, no.~2, p. 154, apr 2024.

\end{thebibliography}



\end{CJK}
\end{document}